\documentclass[11pt]{article}
\usepackage{mathrsfs}
\usepackage{latexsym,lineno}
\usepackage{epsfig}
\usepackage{amsmath}\usepackage{fleqn}\usepackage{verbatim}\usepackage{epsf}
\usepackage{amsthm}\usepackage{graphicx, float}\usepackage{graphicx} \usepackage{tikz}
\usepackage{pgflibraryarrows}
\usepackage{pgflibrarysnakes}
\usepackage{amsfonts}\usepackage{amssymb}\usepackage{graphpap}
\usepackage{epic}\usepackage{curves}

\newtheorem{theorem}{Theorem}[section]
\newtheorem{proposition}[theorem]{Proposition}
\newtheorem{lemma}[theorem]{Lemma}

\newtheorem{ob}{Observation}[section]

\numberwithin{equation}{section}

\title{\bf  Extremal hypergraphs for matching number and domination number\thanks {Research was partially supported by  the National Nature Science Foundation of China (Nos. 11571222, 11471210)}}
\author {Erfang Shan$^{1,2}$, \,Yanxia Dong$^{1}$, \,  Liying Kang$^{1}$\thanks{\em Corresponding authors. Email address: lykang@shu.edu.cn (L. Kang)}, \, Shan Li$^{1}$\\
{\small $^{1}$Department of Mathematics, Shanghai University,
Shanghai 200444, P.R. China}\\
{\small$^{2}$School of Management, Shanghai University,
Shanghai 200444, P.R. China}}

\date{}
\begin{document}

\maketitle

\begin{abstract}
A matching in a hypergraph $\mathcal{H}$ is a set of pairwise disjoint hyperedges. The matching number $\nu(\mathcal{H})$ of $\mathcal{H}$ is the size of a maximum matching in $\mathcal{H}$. A subset $D$ of vertices of $\mathcal{H}$ is  a dominating set of $\mathcal{H}$ if for every $v\in V\setminus D$ there exists $u\in D$ such that $u$ and $v$ lie in an hyperedge of $\mathcal{H}$. The cardinality of a minimum dominating set of $\mathcal{H}$ is  the domination number of $\mathcal{H}$, denoted by $\gamma(\mathcal{H})$. It was proved that $\gamma(\mathcal{H})\leq (r-1)\nu(\mathcal{H})$ for $r$-uniform hypergraphs and the  2-uniform hypergraphs (graphs) achieving equality $\gamma(\mathcal{H})=\nu(\mathcal{H})$ have been characterized. In this paper we  generalize  the inequality $\gamma(\mathcal{H})\leq (r-1)\nu(\mathcal{H})$ to arbitrary  hypergraph  of rank $r$ and we completely characterize  the extremal hypergraphs $\mathcal{H}$  of rank $3$ achieving equality $\gamma(\mathcal{H})=(r-1)\nu(\mathcal{H})$.

\bigskip

\noindent{\bf AMS (2000) subject classification:}  05C65, 05C69, 05C70

\end{abstract}

\section{Introduction}

Hypergraphs are a natural generalization of undirected graphs in which ``edges" may consist of more than 2 vertices.
More precisely, a {\em (finite) hypergraph} $\mathcal{H}=(V, E)$ consists of a (finite) set $V$ and
a collection $E$ of non-empty subsets of $V$. The elements of $V$ are called {\em vertices} and
the elements of $E$ are called {\em hyperedges}, or simply {\em edges} of the
hypergraph. If there is a risk of confusion we will denote the vertex set and the edge set of a hypergraph $\mathcal{H}$ explicitly
by $V(\mathcal{H})$ and $E(\mathcal{H})$, respectively. A hypergraph $\mathcal{H}=(V, E)$  is {\em simple} if no edge is contained in any other edge and $|e|\ge 2$ for all $e\in E$.
 Throughout this paper, we only consider simple hypergraphs.

The {\em rank} of a hypergraph $\mathcal{H}$ is $r(\mathcal{H})={\rm max}_{e\in E} |e|$. An $r$-{\em edge} in $\mathcal{H}$ is an edge of size $r$. The hypergraph $\mathcal{H}$ is said to be $r$-{\em uniform} if every edge of $\mathcal{H}$ is an $r$-edge.
Every (simple) graph is a
2-uniform hypergraph. Thus graphs are special hypergraphs.


A {\em matching} in a hypergraph $\mathcal{H}$ is a set of disjoint hyperedges. The {\em matching number}, denoted by $\nu(\mathcal{H})$, of a hypergraph $\mathcal{H}$ is the size of a maximum
matching in $\mathcal{H}$.

 A subset $D$ of $V(\mathcal{H})$ is called a {\em dominating set} of $\mathcal{H}$ if for every $v\in V(\mathcal{H})\setminus D$ there exists $u\in D$ such that $u$ and $v$ lie in  an hyperedge of $\mathcal{H}$. The minimum cardinality of a dominating set of  $\mathcal{H}$ is called its {\em domination number}, denoted by $\gamma(\mathcal{H})$. Dominating sets are important objects in communication networks, as they represent parts from
which the entire network can be reached directly. Messages can then be transmitted from sources to destinations via such a ``backbone" with suitably chosen links.
The literature on domination has been surveyed
and detailed in \cite{hhs1,hhs2, hy}.
Domination in hypergraphs was introduced  by Acharya \cite{Acharya1} and studied further in \cite{Acharya2,Bujt,Henning1,Jose}.

A subset $T$ of vertices in a hypergraph $\mathcal{H}$ is a {\em transversal} (also called {\em cover}
or {\em hitting set} in many papers) if $T$ has a nonempty intersection with every edge of
$\mathcal{H}$.
The {\em transversal number}, $\tau(\mathcal{H})$, of $\mathcal{H}$ is  the minimum size  of a transversal of $\mathcal{H}$. Transversals in hypergraphs are extensively studied in the literature (see, for example, \cite{Alon,Chv,Dorf,Henning2,Henning3,Lai}).

By definition,   it is easy to see that any transversal of  a hypergraph $\mathcal{H}$ without isolated vertex  is  a dominating set of $\mathcal{H}$ and it must meet all edges of a maximum matching of  $\mathcal{H}$. Furthermore, note that the union of the edges of a maximal matching in $\mathcal{H}$ obviously forms a transversal.  Hence the transversal number of $\mathcal{H}$  is at most $r$ times
its matching number. We state these relationships among the transversal number, the domination number and the matching
number in hypergraphs as an observation.
\begin{ob}
For a hypergraph $\mathcal{H}$ of rank $r$ without isolated vertex, $\nu(\mathcal{H})\le \tau(\mathcal{H})$, $\gamma(\mathcal{H})\leq \tau(\mathcal{H})$,  and $\tau(\mathcal{H})\leq r\nu(\mathcal{H})$.
\end{ob}
Arumugam et al. \cite{Aru} investigated the hypergraphs satisfying $\gamma(\mathcal{H})=\tau(\mathcal{H})$, and  proved that their
recognition problem is NP-hard already on the class of linear hypergraphs of rank 3.
 A long-standing conjecture,
known as
Ryser's conjecture,  asserts that $\tau(\mathcal{H})\leq (k-1)\nu(\mathcal{H})$ for  an $k$-partite  hypergraph $\mathcal{H}$
(see, e.g. \cite{RON,has}). The conjecture turns to be notoriously difficult and remains open for
$k\ge 4$. The relationship between the parameters $\tau(\mathcal{H})$ and $\nu(\mathcal{H})$ in hypergraphs has also been  studied in \cite{Coc,fu,Henning4}.

In particular, if a hypergraph  is $2$-uniform, that is, it is a (simple) graph,
then the following inequality chain is well-known.
\begin{theorem}[\cite{hhs1}] \label{thm01}
If $G$ is a graph  without isolated vertex, then $\gamma(G)\leq \nu(G)\leq \tau(G)$.
\end{theorem}

We observed in  \cite{kld} that the above inequality chain  does not hold in hypergraphs and the difference $\gamma(\mathcal{H})-\nu(\mathcal{H})$  may be arbitrarily large for hypergraphs $\mathcal{H}$ with rank $r\ge 3$.  However, we can extend the inequality $\gamma(G)\leq \nu(G)$ for graphs
 to uniform hypergraphs in \cite{kld}  as follows: for an $r$-uniform hypergraph $\mathcal{H}$ with no isolated vertex, $\gamma(\mathcal{H})\leq (r-1)\nu(\mathcal{H})$,
and this bound is sharp. This paper further observes the the inequality still holds  for  arbitrary hypergraphs of rank $r$.

\begin{theorem}\label{thm02}
If $\mathcal{H}$ is a hypergraph of rank $r$ $(\ge 2)$ without isolated vertex, then $\gamma(\mathcal{H})\leq (r-1)\nu(\mathcal{H})$.
\end{theorem}

In general, a constructive characterization of extremal hypergraphs of rank $r$
achieving equality in Theorem \ref{thm02} seems difficult to obtain. For $2$-uniform hypergraphs, i.e., graphs, Kano et al. \cite{Kano} gave
a complete characterization for extremal graphs with the equality by giving a characterization of star-uniform graphs and showing that  a graph $G$ is star-uniform if and only if $\gamma(G)=\nu(G)$.
In Section 4, we will provide a complete characterization of extremal hypergraphs of rank $3$ with equality $\gamma(\mathcal{H})=(r-1)\nu(\mathcal{H})$.


\section{Terminology and notation}

Let us introduce more definitions and notations. Let $\mathcal{H}=(V(\mathcal{H}), E(\mathcal{H}))$ be  a  hypergraph.
Two vertices $u$ and $v$ are {\em adjacent} in $\mathcal{H}$ if there is an edge $e\in E(\mathcal{H})$ such that $u, v\in e$.
A vertex $v$ and an edge $e$ of $\mathcal{H}$ are {\em incident} if $v\in e$.
The {\em degree}  of a vertex $v\in V(\mathcal{H})$,
denoted by $d_{\mathcal{H}}(v)$ or simply by $d(v)$
if $\mathcal{H}$ is clear from the context, is the number of edges incident to $v$.  A vertex of degree zero is called
an {\em isolated vertex}.
 A vertex of degree $k$ is called a {\em degree-$k$ vertex}.
Let $r,n$ be integers, $1\le r\le n$. We define the $r$-{\em uniform complete hypergraph} on
 $n$ vertices (or the $r$-{\em complete hypergraph}) to be a hypergraph, denoted $K_n^r$, consisting of
all the $r$-subsets of a set of cardinality $n$.

A {\em partial hypergraph} $\mathcal{H}'=(V', E')$ of a hypergraph $\mathcal{H}=(V, E)$, denoted by $\mathcal{H}'\subseteq \mathcal{H}$, is a hypergraph such
that $V'\subseteq V$ and $E'\subseteq E$. In the class of graphs, partial hypergraphs are called {\em subgraphs}.
In particular, if $V'=V$, $\mathcal{H}'$ is called a {\em spanning partial hypergraph} of $\mathcal{H}$.
The partial hypergraph
$\mathcal{H}'=(V', E')$ is {\em induced} if $E'= \{e\in E \mid e\subseteq V'\}$. Induced hypergraphs will be denoted by $\mathcal{H}[V']$. The {\em star} $\mathcal{H}(x)$ centered in $x$ is the family of hyperedges containing $x$.

Two vertices $u$ and $v$ in $\mathcal{H}$ are
{\em connected} if there is a sequence $u=v_{0},v_{1},\ldots,v_{k}=v$ of vertices of $\mathcal{H}$
in which $v_{i-1}$ is adjacent to $v_{i}$ for $i=1,2,\ldots,k$.
A {\em connected hypergraph} is a hypergraph
in which each pair of vertices are connected. A maximal connected partial hypergraph of $\mathcal{H}$ is a
{\em connected component} of $\mathcal{H}$.
Thus, no edge in $\mathcal{H}$ contains vertices from different components.

\section{Proof of Theorem \ref{thm02}}

 A direct proof  of Theorem \ref{thm02} is short.
For the sake of the characterization in next section,
we need to include the  proof.

\noindent{\bf Proof of Theorem \ref{thm02}}\,
We may assume that $\mathcal{H}$ is connected (otherwise we treat each connected component separately).
Let $H^*$ be a hypergraph obtained from $H$ by successively deleting edges of $H$ that do not contain any vertices of degree 1 in the resulting hypergraph at each stage. It is clear that $r(\mathcal{H}^*)\leq r(\mathcal{H})=r$.  When $H$ is transformed to $H^*$, note that isolated vertices cannot arise and the domination number cannot decrease, and every edge of $H^*$ contains at least one degree-1 vertex, so
 every dominating set of $H^*$ is  a transversal of $H^*$. Hence $\tau(H^*)=\gamma(H^*)\ge \gamma(H)$.

 Let
$M=\{e_1, e_2, \ldots, e_l\}$
be a maximum matching of $\mathcal{H}^{*}$. Then $\nu(\mathcal{H}^{*})=|M|$. Let $e_i'$ be the set of vertices obtained from $e_i$ by deleting the vertices of degree-$1$ in $\mathcal{H}^{*}$. We claim that $D=\cup_{i=1}^l e_i'$ is a dominating set of $\mathcal{H}^{*}$. Indeed, for any vertex $x\in V(\mathcal{H}^{*})\setminus D$, $\mathcal{H}^{*}$ has a hyperedge $e$ containing $x$. By the maximality of $M$, $e$ must intersect $V(M)$. This implies that $e\cap D\not=\emptyset$. Thus $D$ is a dominating set of $\mathcal{H}^{*}$.  Hence $\gamma(\mathcal{H}^{*})\leq |D|\leq (r-1)|M|=(r-1)\nu(\mathcal{H}^{*})$.
On the other hand, note that every maximum matching of $\mathcal{H}^{*}$ is  a matching of $\mathcal{H}$, so $\nu(\mathcal{H}^{*})\le \nu(\mathcal{H})$. Therefore, $\gamma(\mathcal{H})\leq \gamma(\mathcal{H}^{*})\leq (r-1)\nu(\mathcal{H}^{*})\leq (r-1)\nu(\mathcal{H})$.
\qed

\section{Hypergraphs $\mathcal{H}$ of rank $3$ with $\gamma(\mathcal{H})=2\nu(\mathcal{H})$}
 In this section
we will give a complete constructive characterization of  hypergraphs $\mathcal{H}$ of rank $3$ satisfying $\gamma(\mathcal{H})=2\nu(\mathcal{H})$.

\subsection{The family $\mathcal{\widehat{H}}_3$ of hypergraphs $\mathcal{H}$ of rank $3$  with $\gamma (\mathcal{H})=2\nu(\mathcal{H})$}
To complete the characterization,
we  construct a family  $\mathcal{\widehat{H}}_3$ of hypergraphs of rank $3$ in which each hypergraph $\mathcal{H}$ satisfies $\gamma (\mathcal{H})=2\nu(\mathcal{H})$.

Let $\mathcal{A}$ be a family of $(2l+1)\times (2l+1)$ upper-triangular matrices,
where $l$, $l\ge 1$, is an arbitrary integer, and for any $A=(a_{ij})\in \mathcal{A}$, $a_{ij}$ is an arbitrary positive integer when $j>i$ and all elements on the main diagonal and elsewhere are zero.
 Let $X=\{x_{1},x_{2}, \ldots, x_{2l+1}\}$ and let $X(A)=\bigcup_{1\leq i<j\leq 2l+1}X_{ij}$, where $X_{ij}=\{x^{p}_{ij}\mid p=1,2, \ldots, a_{ij},  a_{ij}\in A\}$,
 for $1\leq i<j\leq 2l+1,  A\in \mathcal{A}$.
 A family of $3$-uniform hypergraphs $\mathcal{H}_3(A)$ is defined as
$$V(\mathcal{H}_3(A))=X\cup X(A),\, \,  E(\mathcal{H}_3(A))=\bigcup_{1\leq i<j\leq 2l+1}E_{ij},$$
where $E_{ij}=\big\{\{x_i,v,x_j\}\mid v\in X_{ij}\big\}$ and $A\in \mathcal{A}$. Set $S=X(A)-\{x_{ij}^1\mid 1\leq i<j\leq 2l+1\}$.
Let $\mathcal{E}_1={X\choose 3}\bigcup{X\choose 2}$ where ${X\choose r}$ denote
the family of all  $r$-subsets of $X$,
$\mathcal{E}_2=\big\{\{x_i,x_j,v\}\mid 1\leq i<j\leq 2l+1, v\in S\big\}$.
Finally, we define the family $\mathcal{\widehat{H}}_3$ of hypergraphs of rank $3$.
For every $\mathcal{H}\in \mathcal{\widehat{H}}_3$, $\mathcal{H}$ is defined as
$$V(\mathcal{H})=X\cup X(A),\, \, E(\mathcal{H})=\big\{E(\mathcal{H}_3(A))\cup \mathcal{E}_1'\cup \mathcal{E}_2'
\mid \mathcal{E}_1'\subseteq \mathcal{E}_1, \mathcal{E}_2'\subseteq \mathcal{E}_2\big\}.$$

Furthermore, we define a hypergraph  $\mathcal{F}$, where $V(\mathcal{F})=\{x_1, x_2, x_3, x_{12}^1, x_{23}^1, x_{13}^1\}$,
$E(\mathcal{F})=\big\{\{x_1,x_{12}^1, x_2\}, \{x_2, x_{23}^1, x_3\}, \{x_1, x_{13}^1, x_3\}, \{x_{12}^1, x_{23}^1, x_{13}^1\}\big\}$.

{\bf Example.} Let
\begin{equation*}
A= \left( {\begin{array}{ccc}
   0 & 2 & 1 \\
   0 & 0 & 2 \\
   0 & 0 & 0 \\
\end{array}} \right).
\end{equation*}
The hypergraph $\mathcal{F}$, $\mathcal{H}_3(A)$ and an example of hypergraphs  in $\mathcal{\widehat{H}}_3$ constructed by $X=\{x_1,x_2,x_3\}$ and $A$ are exhibited in Fig. 1.

\begin{figure}[htbp]
  \centering
  \includegraphics[height=5.5cm]{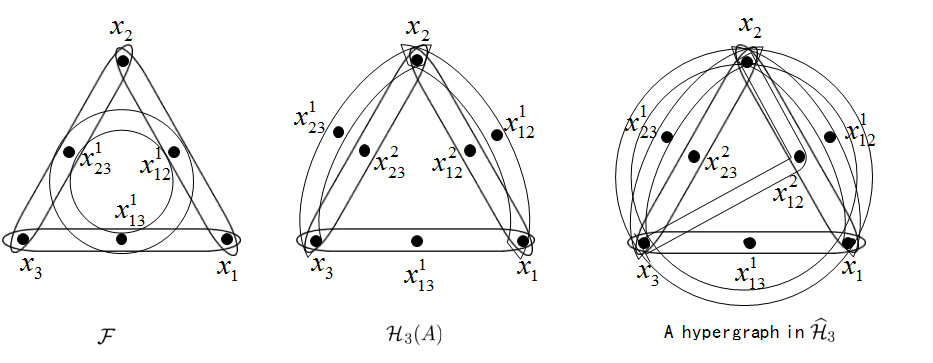}
  \caption{The hypergraph $\mathcal{F}$, $\mathcal{H}_3(A)$ and an example of hypergraphs  in $\mathcal{\widehat{H}}_3$} \label{Fig1}
\end{figure}

For $\mathcal{H}=\mathcal{F}$ or  $\mathcal{H}\in \mathcal{\widehat{H}}_3$, by the above construction, $X-\{x_i\}$
is a minimum dominating set of $\mathcal{H}$ for any $x_i\in X$, and $\big\{\{x_1, x_{12}^1, x_2\}, \{x_3, x_{34}^1, x_4\},\ldots,
\{x_{2l-1}, x_{(2l-1)2l}^1, x_{2l}\}\big\}$ is a maximum matching of  $\mathcal{H}$.
 So we  immediately
have the following proposition.

\begin{proposition}\label{pro3.0}
For $\mathcal{H}=\mathcal{F}$ or  $\mathcal{H}\in \mathcal{\widehat{H}}_3$, $\nu(\mathcal{H})=l,
\gamma(\mathcal{H})=2l$ and $\gamma(\mathcal{H})=2\nu(\mathcal{H})$.
\end{proposition}

Furthermore,  we have the following property.

\begin{proposition}\label{pro3.1}
Let $\mathcal{H}=\mathcal{F}$ or  $\mathcal{H}\in \mathcal{\widehat{H}}_3$.
For any two vertices $x,y\in V(\mathcal{H})$,  if $x$ and $y$ satisfy  one of the following conditions:

{\rm (i)} $x,y\in X(A)$,

{\rm (ii)} $|X\cap \{x,y\}|=1$ and $|X|\geq 5$,

{\rm (iii)} $|X\cap \{x,y\}|=1$, $|X|=3$ and $x$ and $y$ are adjacent,\\
then there  exists a maximum matching $M$ of $\mathcal{H}$ such that $x,y\not\in V(M)$.
\end{proposition}
\noindent{\bf Proof.}
If $\mathcal{H}=\mathcal{F}$, clearly the assertion holds. So we may assume that $\mathcal{H}\in \mathcal{\widehat{H}}_3$. We show the assertion by consider the conditions separately.

Suppose that $x,y\in X(A)$.
Without loss of generality, we may assume that $x$ and $x_1$ are adjacent.
Clearly,
$M=\big\{\{x_2,x_{23}^1,x_3\}, \{x_4,x_{45}^1,x_5\},\ldots, \{x_{2l},x_{2l(2l+1)}^1, x_{2l+1}\}\}$ is a maximum matching of
$\mathcal{H}$, where $l=\nu(\mathcal{H})$.
If $y\not\in V(M)$, we are done.
Otherwise, there exists an integer $i$, $1\le i\le l$, such that $y=x_{2i(2i+1)}^1$,
that is, $y$ lies in the edge $\{x_{2i},x_{2i(2i+1)}^1, x_{2i+1}\}$
of $M$.
If $|X|=3,$ then $\nu(\mathcal{H})=l=1$, hence
$\{x_1, x_{12}^1, x_2\}$ or $\{x_1, x_{13}^1, x_3\}$ is a maximum matching of $\mathcal{H}$ satisfying $x,y\not\in V(M)$.
If  $|X|\geq 5$,
there exists an integer $j\neq i, 1\leq j\leq  l$ such that $\{x_{2j},x_{2j(2j+1)}^1, x_{2j+1}\}$ and $\{x_{2i},x_{2i(2i+1)}^1, x_{2i+1}\}$
are two distinct edges of $M$.  Without loss of generality, let $i<j$.
We set
\begin{align*}
M'=&M\cup \big\{\{x_{2i},x_{2i(2j+1)}^1,x_{2j+1}\}, \{x_{2i+1},x_{(2i+1)(2j)}^1, x_{2j}\}\}\\
&-\big\{\{x_{2j},x_{2j(2j+1)}^1, x_{2j+1}\}, \{x_{2i},x_{2i(2i+1)}^1, x_{2i+1}\}\big\}.
\end{align*}
Then $M'$ is a maximum matching of $\mathcal{H}$ satisfying $x,y\not\in V(M)$.

Suppose that $|X\cap \{x,y\}|=1$, and either $|X|\geq 5$ or $|X|=3$ and $x$ and $y$ are adjacent. Without loss of generality, we may assume that $x=x_1$.  In this case,
its proof is similar to the above proof, we can prove that the assertion is true.
\qed

\subsection{Extremal hypergraphs  of rank $3$ achieving  $\gamma(\mathcal{H})=2\nu(\mathcal{H})$}

We  define
$$\mathfrak{H}_3=\{\mathcal{H}\,|\, \mathcal{H} \, \, \mbox{has rank $3$ and satsifies that $\gamma (\mathcal{H})=2\nu(\mathcal{H})$} \}.$$

The following lemma lists various basic properties on hypergraphs in $\mathfrak{H}_3$.

\begin{lemma}\label{lem3.1}
For every hypergraph $\mathcal{H}\in \mathfrak{H}_{3}$, there exists a
spanning partial hypergraph $\mathcal{H}^{*}$ of $\mathcal{H}$ satisfying the following properties.

{\rm (i)} $\mathcal{H}^*\in \mathfrak{H}_{3}$, $\gamma(\mathcal{H}^{*})=\gamma(\mathcal{H})$ and $\nu(\mathcal{H}^*)=\nu(\mathcal{H})$.

{\rm (ii)} $\mathcal{H}^{*}$ is a $3$-uniform hypergraph.

{\rm (iii)} For each connected component $\mathcal{H}^*_{i}$ of $\mathcal{H}^{*}$, $\mathcal{H}^*_i\in \mathfrak{H}_{3}$.

{\rm (iv)} Every edge in $\mathcal{H}^{*}$ contains exactly one degree-$1$ vertex.

\end{lemma}

\noindent{\bf Proof.}
As  described in the proof of Theorem \ref{thm02}, there exists a spanning partial hypergraph  $\mathcal{H}^{*}$ (not necessarily connected) of $\mathcal{H}$
 such that every edge of $\mathcal{H}^{*}$ contains at least one degree-$1$ vertex.
Let $\mathcal{H}^*_{1}, \mathcal{H}^*_{2}, \ldots, \mathcal{H}^*_{k}$ be the connected  components of $\mathcal{H}^*$. Next we show that  $\mathcal{H}^{*}$
is the spanning partial hypergraph satisfying the properties (i)-(iv).

As noted in the proof of Theorem \ref{thm02}, $\gamma(\mathcal{H})\leq \gamma(\mathcal{H}^{*})$ and $\nu(\mathcal{H}^*)\leq \nu(\mathcal{H})$, and $\gamma(\mathcal{H}^{*})\leq (r(\mathcal{H}^{*})-1)\nu(\mathcal{H}^{*})\leq 2 \nu(\mathcal{H}^{*})$.
Note that $2\leq r(\mathcal{H}^{*})\leq3$. Since $\mathcal{H}\in \mathfrak{H}_3$, we have
$$2 \nu(\mathcal{H})=\gamma(\mathcal{H})\leq \gamma(\mathcal{H}^*)\leq (r(\mathcal{H}^{*})-1)\nu(\mathcal{H}^{*})\leq 2\nu(\mathcal{H}^*)\leq 2 \nu(\mathcal{H}),$$
and hence $\nu(\mathcal{H}^*)= \nu(\mathcal{H})$ and $\gamma(\mathcal{H}^*)=\gamma(\mathcal{H})$.
Further it is obtained that $\gamma(\mathcal{H}^*)=2\nu(\mathcal{H}^*)$ and $r(\mathcal{H}^{*})=3$. Therefore, $\mathcal{H}^{*}\in \mathfrak{H}_{3}$, and part (i) follows.

Next we show that $\mathcal{H}^{*}$ is a $3$-uniform hypergraph. Since $r(\mathcal{H}^{*})=3$, suppose that
there exists an edge $e\in E(\mathcal{H}^{*})$ such that $|e|=2$.
Let $M$ be a maximum matching of $\mathcal{H}^*$ and $U$ the set of degree-$1$ vertices of $\mathcal{H}^*$.
Then $\nu(\mathcal{H}^*)=|M|$ and $|V(M)|\leq 3\nu(\mathcal{H}^*)$.
As we have seen in the proof of Theorem \ref{thm02}, $V(M)\backslash U$ is a dominating set of $\mathcal{H}^*$,
i.e., $\gamma(\mathcal{H}^*)\leq|V(M)\backslash U|\leq2|V(M)|/3\leq 2\nu(\mathcal{H}^*)$.
So it follows from $\mathcal{H}^{*}\in \mathfrak{H}_{3}$ that
 each edge of $M$ has $3$ of the size. Thus $e \not\in M$. Since $M$ is a maximum matching of $\mathcal{H}^*$, we have $e\cap M\neq \emptyset$.
 By the proof of Theorem \ref{thm02}, $e$ contains exactly one degree-1 vertex. This implies that $e$ meets $M$ precisely in one edge, say $f$.
 Hence $M'=e\cup (M\backslash f)$ is also a maximum matching of $\mathcal{H}^*$.
Thus, $V(M')\backslash U$ would be a dominating set of $\mathcal{H}^*$, i.e.,
$\gamma(\mathcal{H}^*)\leq|V(M')\backslash U|\leq 2\nu(\mathcal{H}^*)-1$, contradicting the fact that $\gamma(\mathcal{H}^*)=2\nu(\mathcal{H}^*)$.
The part (ii) follows.

Since each connected component $\mathcal{H}^*_{i}$ of $\mathcal{H}^*$ is also a $3$-uniform hypergraph,
$\gamma(\mathcal{H}^*_{i})\leq 2 \nu(\mathcal{H}^*_{i})$ by Theorem \ref{thm02}.
Hence by the equality $\gamma(\mathcal{H}^*)= 2 \nu(\mathcal{H}^*)$, we have
$$\sum_{i=1}^{k}\gamma(\mathcal{H}^*_{i})=\gamma(\mathcal{H}^*)=2 \nu(\mathcal{H}^*)=2\sum_{i=1}^{k}\nu(\mathcal{H}^*_{i}).$$
 Combining this with $\gamma(\mathcal{H}^*_{i})\leq 2 \nu(\mathcal{H}^*_{i})$,  we  obtain $\gamma(\mathcal{H}^*_{i})=2\nu(\mathcal{H}^*_{i})$, that is, $\mathcal{H}^*_{i} \in \mathfrak{H}_{3}$ for each $1\le i\le k$. Then $\mathcal{H}^*_{i}$ satisfies the properties (iii).

Now we show that every edge in $\mathcal{H}^{*}$ contains exactly one degree-$1$ vertex. It suffices to prove the assertion for each connected component
of $\mathcal{H}^{*}$. Let $\mathcal{H}^*_{i}$ be an arbitrary connected component of $\mathcal{H}^{*}$.

Suppose, to the contrary, that $\mathcal{H}^*_i$ has an edge $e$ that  contains at least two degree-$1$
vertices. Let $x_{1}, x_{2}\in e$ and $d_{\mathcal{H}^*_i}(x_1)=d_{\mathcal{H}^*_i}(x_2)=1$.
If each vertex of  $e$ is a  degree-$1$ vertex,
 then  $\mathcal{H}^*_i$  consists of the single edge $e$, and hence
 $\nu(\mathcal{H}^*_i)=\gamma(\mathcal{H}^*_i)=1$, contradicting the fact that $\gamma(\mathcal{H}^*_i)=2\nu(\mathcal{H}^*_i)$.
 We may
assume, therefore, that $e$ contains a vertex, say $x$, such that $d_{\mathcal{H}^*_i}(x)\ge 2$. Thus $e=\{x, x_{1}, x_{2}\}$.

Let $M_i=\{f_{1}, f_{2}, \ldots, f_{\nu(\mathcal{H}^*_i)}\}$ be a maximum matching of the connected component $\mathcal{H}^*_i$ of $\mathcal{H}^*$.
Let $I=\{y\in \mathcal{H}^*_i \mid \  d_{\mathcal{H}^*_i}(y)=1\}$. It is easy to see that $D=V(M_i)\setminus I$ is a dominating set of $\mathcal{H}^*_i$. This implies  that each edge in $M_i$
 contains exactly one degree-1 vertex of $\mathcal{H}^*_i$, for otherwise $\gamma(\mathcal{H}^*_i)\leq |D|\leq 2\nu(\mathcal{H}^*_i)-1$, contradicting the fact that $\gamma(\mathcal{H}^*_i)=2\nu(\mathcal{H}^*_i)$.
Therefore, $e\in E(\mathcal{H}^*_i)\setminus M_{i}$.
By the maximality of $M_{i}$, there exists an edge $f\in M_i$ such that $f\cap e\not=\emptyset$. Since $x$ is the unique vertex of  degree at least 2 in $e$, $x\in f\cap e$. Let $f=\{u, v, x\}$. As noted above, $f$ contains exactly one degree-1 vertex. Without loss of generality, let $d_{\mathcal{H}^*_i}(v)\ge 2$.
Then $|E(\mathcal{H}_i^*(v))|\ge 2$. Let $g\in E(\mathcal{H}^*_i(v))$ and $g\neq f$. One sees that $g\neq e$, and either
 $g\cap f=\{v\}$ or $g\cap f=\{v,x\}$.
 If  $g\cap f=\{v\}$, then $g\cap e=\emptyset$. This implies that $g\cap V(M_{i}\setminus \{f\})\neq\emptyset$, for otherwise we obtain a matching $(M_{i}\setminus \{f\})\cup\{g, e\}$, contradicting the maximality of $M_{i}$. Therefore, $g\cap f=\{v\}$ and $g\cap V(M_{i}\setminus \{f\})\neq\emptyset$ or $g\cap f=\{v,x\}$.
 Thus $D=V(M_i)\setminus (I\cup\{v\})$ is a dominating set of $\mathcal{H}_{i}^*$. But then $\gamma(\mathcal{H}^*_i)\leq2\nu(\mathcal{H}^*_i)-1$, a contradiction.
The part (iv) follows. \qed

In what follows, let $\mathcal{H}^*$ be the spanning partial hypergraph of $\mathcal{H}$ satisfying the properties in Lemma \ref{lem3.1},
$\mathcal{H}^*_{1}, \mathcal{H}^*_{2}, \ldots, \mathcal{H}^*_{k}$  all connected  components of $\mathcal{H}^*$, and $\mathcal{H}_i=\mathcal{H}[V(\mathcal{H}^*_i)]$ for $1\le i\le k$.

For each connected component $\mathcal{H}^*_i$ of $\mathcal{H}^*$, we define an {\em edge-contracting hypergraph} $G_i$ of $\mathcal{H}^*_i$ as follows.
Let $G_i$ be the hypergraph obtained from $\mathcal{H}^*_i$ by deleting all the degree-$1$ vertices of $\mathcal{H}^*_i$
so that each 3-edge of $\mathcal{H}^*_i$ is contracted into a 2-edge, and deleting all multi-edges, if any, from the resulting hypergraph.
Since $\mathcal{H}^*_i$ is a connected $3$-uniform hypergraph,
the resulting hypergraph $G_i$ is a connected graph.

For the edge-contracting hypergraph $G_i$ of  connected component $\mathcal{H}^*_i$ of $\mathcal{H}^*$, we show that
$G_i$ is a complete graph of odd order.
\begin{lemma}\label{lem3.2}
For  every connected component $\mathcal{H}^*_i$ of $\mathcal{H}^*$,  the edge-contracting hypergraph  $G_i$ of $\mathcal{H}^*_i$ is a complete graph on $2\nu(\mathcal{H}_i^*)+1$ vertices.
\end{lemma}
\noindent{\bf Proof.}
Let  $M_i$ be  a maximum matching of $\mathcal{H}^*_i$ and $I$  the set of degree-$1$ vertices of $\mathcal{H}^*_i$.
Then $\nu(\mathcal{H}^*_i)=|M_i|$.
Let $Y=V(\mathcal{H}^*_i)\setminus I$, $Y_1=V(M_i)\cap Y$. Clearly, $Y=V(G_i)$ and $Y_1=V(M_i)\setminus I$.
As we have seen in the proof of Theorem \ref{thm02}, $Y_1$ $(=V(M_i)\setminus I)$ is a minimum dominating set of $\mathcal{H}^*_i$, i.e., $\gamma(\mathcal{H}^*_i)=|Y_1|=2\nu(\mathcal{H}^*_i)=2|V(M_i)|/3$.
Furthermore, we have the following claims.

\noindent
{\bf Claim 1.}
$|Y\setminus Y_1|=1$ and every vertex of $Y_1$ is adjacent to the vertex of $Y\setminus Y_1$.

\noindent{\bf Proof of Claim 1.}
We establish the claim by contradiction.
Suppose first that $Y\setminus Y_1=\emptyset$.
Choose an arbitrary vertex $v$ of $Y_1$. We show that  $Y_1\setminus \{v\}$ is a dominating set of $\mathcal{H}^*_i$.
Indeed, for each vertex $x$ in $V(\mathcal{H}^*_i)\setminus V(M_i)$, there is an edge $e\in E(\mathcal{H}^*_i)\setminus M_i$
such that $e$ contains $x$.
$Y\setminus Y_1=\emptyset$ implies that $Y=Y_1$ and $V(\mathcal{H}^*_i)\setminus V(M_i)\subseteq I$. By Lemma \ref{lem3.1},
every edge in $\mathcal{H}^{*}_i$ contains exactly one degree-$1$ vertex, so $x$ is the unique degree-1 vertex in $e$.
It follows that $|e\cap (V(M_i)\setminus I)|=|e\cap Y_1|\ge 2$. Hence $|e\cap (Y_1\setminus \{v\}|\ge 1$, that is, $x$ is dominated by $Y_1\setminus \{v\}$.
Therefore, $Y_1\setminus \{v\}$ is a dominating set of $\mathcal{H}^*_i$,  contradicting the fact that $\gamma(\mathcal{H}^*_i)=|Y_1|$.

Suppose that $|Y\setminus Y_1|\geq2$. Let $v_{1}$ and  $v_{2}$ be any two distinct vertices in  $Y\setminus Y_1$, and let $e_1$ and  $e_2$ be two edges of $\mathcal{H}^*_i$ containing $v_1$ and $v_2$, respectively. Clearly, $e_1, e_2\not\in M_i$.
Since $M_i$  is a maximum matching  of $\mathcal{H}^*_i$, $e_j\cap Y_1\neq \emptyset$ for $j=1,2$. By Lemma \ref{lem3.1}, every edge in $\mathcal{H}^{*}_i$ contains exactly one degree-$1$ vertex, so $e_1\neq e_2$, $v_1\not\in e_2$ and $v_2\not\in e_1$.
By the maximality of $M_i$, there exist two edges $f_1$ and $f_2$ in $M_i$ such that $f_1\cap e_1\neq\emptyset$ and $f_2\cap e_2\neq\emptyset$.
Let $w_1\in f_1\cap e_1$ and $w_2\in f_2\cap e_2$. By Lemma \ref{lem3.1}, each one of $\{e_1,e_2,f_1, f_2\}$ has exactly one degree-1 vertex.
 Let $e_1=\{u_1, v_1, w_1\}, e_2=\{u_2, v_2, w_2\}$, $f_1=\{x_1, y_1, w_1\}$, and  $f_2=\{x_2, y_2,w_2\}$, where $d_{\mathcal{H}^*_i}(u_1)=d_{\mathcal{H}^*_i}(u_2)=d_{\mathcal{H}^*_i}(x_1)=d_{\mathcal{H}^*_i}(x_2)=1$.
 Then $y_1,y_2\in Y_1$.
It is clear that both $y_1\not=w_2$ and $y_2\neq w_1$, since
 otherwise  $M_i\cup \{e_1,e_2\}\setminus \{f_1\}$ or $M_i\cup \{e_1,e_2\}\setminus \{f_2\}$
would be a matching of $\mathcal{H}^*_i $, contradicting the maximality of  $M_i$.
 We claim that there exists an edge $g_1\in \mathcal{H}_i^*(y_1)$ such that $ (V(M_i)\setminus \{y_1\})\cap g_1=\emptyset$.
  Otherwise, $D=V(M_i)\setminus (I\cup \{y_1\})$ would be a dominating set of $\mathcal{H}^*_i$, contradicting
 the minimality of $Y_1$ ($=V(M_i)\setminus I$).
Furthermore, we have $v_1\in g_1$ for otherwise $(M_i\setminus \{f_1\})\cup \{g_1,e_1\}$ is a matching of  $\mathcal{H}^*_i$,
contradicting the maximality of $M_i$.
If $w_1=w_2$,  similarly we have $v_2\in g_1$, then $g_1$ contains no degree-$1$ vertex, contradicting Lemma \ref{lem3.1}.
So we may assume that $w_1\not= w_2$. Similarly, we see that
there exists an edge $g_2\in \mathcal{H}_i^*(y_2)$ such that $(V(M_{i})\setminus \{y_2\})\cap g_2=\emptyset$ and $v_2\in g_2$.
By the connectivity of $\mathcal{H}^*_i$, there exists an  $M_i$-augmenting path from $v_1$ to $v_2$, again contradicting  the maximality of $M_i$.
Therefore,  $|Y\setminus Y_1|=1$.

 Let $Y\setminus Y_1=\{x\}$. We show that every vertex in $Y_1$ is adjacent to the vertex $x$. Suppose not, let $v\in Y_1$ such that $v$ is not adjacent to $x$.
 By Lemma \ref{lem3.1} (iii), for each vertex $y\in V(\mathcal{H}_i^*)\setminus V(M_i)$ and each edge $e$ containing $y$,
 we have $e\cap (Y_1\setminus \{v\})\neq\emptyset$.
 Then $Y_1\setminus \{v\}$ is a dominating set of $\mathcal{H}^*_i$ with cardinality  smaller than $|Y_1|$, which is a contradiction. This completes the proof of Claim 1. ~$\Box$

\noindent
{\bf Claim 2.}
Any two vertices  in $Y_1$ are adjacent.

\noindent{\bf Proof of Claim 2.}
Suppose to the contrary that there exist two vertices $u, v\in Y_1$  such that $u, v$ are not adjacent.
 By Lemma \ref{lem3.1} (iii),  we deduce that $D=(Y_1\setminus\{u, v\})\cup \{x\}$ is a dominating set of $\mathcal{H}_i^*$, contradicting
 the minimality of $|Y_1|$.  ~$\Box$

 By Claims 1, 2,  it is easily seen  that $G_i=K_{2\nu(\mathcal{H}^*_i)+1}$.
\qed

\begin{lemma}\label{lem3.3}
For  every connected component $\mathcal{H}^*_i$ of $\mathcal{H}^*$,  we have
$\gamma(\mathcal{H}_i)=2\nu(\mathcal{H}_i)$, $\nu(\mathcal{H}_i)= \nu(\mathcal{H}_i^*)$, $\gamma(\mathcal{H}_i)=\gamma(\mathcal{H}_i^*)$,
$\gamma(\mathcal{H})=\sum_{i=1}^{k}\gamma(\mathcal{H}_i)$ and $\nu(\mathcal{H})=\sum_{i=1}^{k}\nu(\mathcal{H}_i)$.
\end{lemma}
\noindent{\bf Proof.} Clearly, $\mathcal{H}_i$ is a hypergraph of rank $3$. By Theorem \ref{thm02}, $\gamma(\mathcal{H}_i)\leq  2\nu(\mathcal{H}_i)$,
so
$$\gamma(\mathcal{H})\le \sum_{i=1}^{k}\gamma(\mathcal{H}_i)\le\sum_{i=1}^{k}2\nu(\mathcal{H}_i)\leq 2\nu(\mathcal{H}).$$
Since $\mathcal{H}\in \mathfrak{H}_3$, $\gamma(\mathcal{H})=2\nu(\mathcal{H})$, and hence
$\gamma(\mathcal{H}_i)=2\nu(\mathcal{H}_i)$ for each $1\le i\le k$. Since   $\gamma(\mathcal{H}_i)\leq \gamma(\mathcal{H}_i^*)=2\nu(\mathcal{H}^*_i)\leq 2\nu(\mathcal{H}_i)$, $\nu(\mathcal{H}_i)= \nu(\mathcal{H}_i^*)$ and $\gamma(\mathcal{H}_i)=\gamma(\mathcal{H}_i^*)$. This completes the proof of the lemma.
\qed

\begin{lemma}\label{thm3.1}
For every $\mathcal{H}\in \mathfrak{H}_3$,  let $\mathcal{H}^*$ be the spanning partial hypergraph of $\mathcal{H}$  described in Lemma \ref{lem3.1}, $\mathcal{H}_1^*,
 \mathcal{H}_2^*, \ldots,  \mathcal{H}_k^*$ be all the connected components of $\mathcal{H}^*$, and let $\mathcal{H}_i=\mathcal{H}[V(\mathcal{H}^*_i)]$ for $1\le i\le k$. Then either
$\mathcal{H}_i=\mathcal{F}$ or   $\mathcal{H}_i\in \mathcal{\widehat{H}}_{3}$ for $1\le i\le k$.
\end{lemma}

\noindent{\bf Proof.}
By Lemma \ref{lem3.2}, the edge-contracting hypergraph $G_i$ of $\mathcal{H}^*_i$ is a complete graph on $2\nu(\mathcal{H}_i^*)+1$ vertices,
so $\mathcal{H}_i^*\in \mathcal{H}_3(A)$. Let $V(G_i)=\{x_{1}, x_{2}, \ldots, x_{2l_i+1}\}$ where $l_i=\nu(\mathcal{H}^*_i)$.
By Lemma \ref{lem3.3}, $\nu(\mathcal{H}_i)=l_i$. For notational simplicity,
we write $X_i$ for $V(G_i)$, let $X_i(A)=\cup_{1\leq i_s<i_t\leq 2l_i+1}X_{i_si_t}$,
where $X_i$ and $X_i(A)$ can be regarded as the sets corresponding to $X$ and $X(A)$, respectively, described in Subsection 3.1.

\noindent
{\bf Claim 3.} If $\mathcal{H}_i\not=\mathcal{F}$, then  $|e\cap X_i|\geq 2$ for each $e\in E(\mathcal{H}_i)\setminus E(\mathcal{H}^*_i)$.

\noindent{\bf Proof of Claim 3.} \ Suppose, to the contrary, that there exists an edge
$e\in E(\mathcal{H}_i)\setminus E(\mathcal{H}^*_i)$ such that $|e\cap X_i|\leq 1$.
We distinguish two cases.

 {\it Case} 1. $|e\cap X_i|=0$. Clearly, $e\subseteq X_i(A)$. Let $e=\{x_{i_1j_1}^{k_1}, x_{i_2j_2}^{k_2}, x_{i_3j_3}^{k_3}\}$.
If $|X_i(A)|=3$, then $|X_i|=3$ and $\nu(\mathcal{H})=1$. So $\mathcal{H}_i=\mathcal{F}$, a contradiction.
So we may assume that  $|X_i(A)|\geq5$. By the construction of $\mathcal{H}_i^*$, it is easy to see that
there exists a maximum matching $M_i$ of $\mathcal{H}^*_i$ such that
the three edges $\{x_{i_1}, x_{i_1j_1}^{k_1}, x_{j_1}\}, \{x_{i_2}, x_{i_2j_2}^{k_2}, x_{j_2}\}$ and  $\{x_{i_3}, x_{i_3j_3}^{k_3}, x_{j_3}\}$
do not belong to $M_i$. This implies that $M_i\cup \{e\}$ is a maximum matching of $\mathcal{H}^*_i\cup \{e\}$. Hence
$\nu(\mathcal{H}^*_i\cup \{e\})=\nu(\mathcal{H}^*_i)+1$.
It follows immediately that $\nu(\mathcal{H}_i)\ge \nu(\mathcal{H}^*_i\cup \{e\})>\nu(\mathcal{H}^*_i)$.
But, by Lemma \ref{lem3.3}, we have $\nu(\mathcal{H}_i)=\nu(\mathcal{H}^*_i)$, a contradiction.

 {\it Case} 2. $|e\cap X_i|=1$. Then $|e\cap X_i(A)|=2$. Without loss of generality, let $e=\{x_{1}, x_{i_1j_1}^{k_1}, x_{i_2j_2}^{k_2}\}$.
 Clearly, $M_i=\big\{\{x_2,x_{23}^1,x_3\}, \{x_4,x_{45}^1,x_5\},\ldots, \{x_{2l_i},x_{2l_i(2l_i+1)}^1,x_{2l_i+1}\}\big\}$ is a maximum matching of
$\mathcal{H}^*_i$  and $X_i\setminus V(M_i)=\{x_{1}\}$.
If $x_{i_1j_1}^{k_1}, x_{i_2j_2}^{k_2}\not\in V(M_i)$, then $M_i\cup \{e\}$ is a matching of $\mathcal{H}_i$, which  contradicts
  Lemma \ref{lem3.3}. Without loss of generality, we may assume that $x_{i_1j_1}^{k_1}=x_{i_1j_1}^1\in V(M_i)$ and   $x_{i_2j_2}^{k_2}=x_{i_2j_2}^1\in V(M_i)$.
 Then $j_1=i_1+1$ and $j_2=i_2+1$.
  Let
   $$M'_i=M_i\cup \big\{\{x_{i_1},x_{i_1j_2}^1,x_{j_2}\}, \{x_{i_2},x_{i_2j_1}^1,x_{j_2}\}\big\}\cup \{e\}\setminus \big\{\{x_{i_1},x_{i_1j_1}^1,x_{j_1}\}, \{x_{i_2},x_{i_2j_2}^1,x_{j_2}\}\big\}.$$
  It is clear that $M'_i$ is a matching of $\mathcal{H}_i$, which contradicts Lemma \ref{lem3.3} again.

Therefore, $|e\cap X_i|\geq 2$. This completes the proof of Claim 3. ~$\Box$

Now we show that $\mathcal{H}_i=\mathcal{F}$ or   $\mathcal{H}_i\in \mathcal{\widehat{H}}_{3}$.
If $\mathcal{H}_i=\mathcal{F}$, then we are done. Otherwise, we show that $\mathcal{H}_i\in \mathcal{\widehat{H}}_{3}$.
Suppose to the contrary that $\mathcal{H}_i\not\in \mathcal{\widehat{H}}_{3}$.
Then there exist two integers $i_s$ and $i_t$, $1\leq i_s<i_t\leq 2l_i+1$ such that $d_{\mathcal{H}_i}(v)\geq 2$ for all $v\in X_{i_si_t}$.
This implies that  $X_i\setminus\{x_{i_s},x_{i_t}\}$ is a dominating set of $\mathcal{H}_i$, $\gamma(\mathcal{H}_i)\le |X_i|-2=
|V(G_i)|-2=2\nu(\mathcal{H}_i)-1$.
But, by Lemma \ref{lem3.3}, $\gamma(\mathcal{H}_i)=2\nu(\mathcal{H}_i)$, a contradiction.
 By Claim 1, we conclude that $\mathcal{H}\in \mathcal{\widehat{H}}_3$.  \qed

As defined in Subsection 4.1, let $A_i\in \mathcal{A}$ and, $X_i$, $X(A_i)$ and $S_i$ corresponding to $A_i$ are the set of vertices,  for $i=1,2,\ldots,  k$.
Write $X=\cup_{i=1}^kX_i$ and let
\begin{eqnarray*}
\mathcal{Z}_1&=&\big\{\{x_{m_1},x_{m_2},x_{n}\}\mid x_{m_1},x_{m_2} \in X_{m}, x_n\in X_{n},m\neq n, 1\le m,n\le k\big\},\\
\mathcal{Z}_2&=&\big\{\{x_{m_1},x_{m_2},y\} \mid x_{m_1},x_{m_2} \in X_{m}, y\in S_{n},m\neq n, 1\le m,n\le k \big\}.
\end{eqnarray*}
We define the family $\mathcal{G}_3$ of 3-uniform hypergraphs as follows: For every $\mathcal{H}\in \mathcal{G}_3$, $\mathcal{H}$ is defined as
\begin{align*}
V(\mathcal{H})&=\bigcup_{i=1}^k(X_i\cup X(A_i)),\\
 E(\mathcal{H})&=\Big\{\bigcup_{i=1}^{k}E(\mathcal{H}_{i})\cup \mathcal{Z}_1'\cup \mathcal{Z}_2'\mid \mathcal{H}_{i}\in \mathcal{\widehat{H}}_{3}
 \cup\{\mathcal{F}\},  \mathcal{Z}_1'\subseteq \mathcal{Z}_1, \mathcal{Z}_2'\subseteq \mathcal{Z}_2\Big\},
\end{align*}
 where $V(\mathcal{H}_{i})=X_{i}\cup X_{i}(A_{i})$.

We are now ready to give a constructive characterization of hypergraphs in $\mathcal{\mathfrak{H}}_3$.
\begin{theorem}\label{thm3.7}
For  a hypergraph $\mathcal{H}$ of rank $3$,  $\gamma(\mathcal{H})=2\nu(\mathcal{H})$ if and only if
$\mathcal{H}\in \mathcal{G}_3$.
\end{theorem}

\noindent{\bf Proof.}
First, suppose that
$\mathcal{H}\in \mathcal{G}_3$. We show that $\mathcal{H}\in \mathfrak{H}_3$.
 By the construction of $\mathcal{G}_3$, $\mathcal{H}$
is a  hypergraph of rank $3$. It remains to  show that $\gamma(\mathcal{H})=2\nu(\mathcal{H})$.
By the definition of $\mathcal{H}$, we have
that
$\gamma(\mathcal{H})\leq \sum_{i=1}^{k}\gamma(\mathcal{H}_{i})$ and
$\nu(\mathcal{H})\geq\sum_{i=1}^{k}\nu(\mathcal{H}_{i})$.
Since $\mathcal{H}_i\in \mathcal{\widehat{H}}_3\cup\{\mathcal{F}\}$, by Proposition \ref{pro3.0}, $\gamma(\mathcal{H}_i)=2\nu(\mathcal{H}_i)$.
Hence,
 $$\gamma(\mathcal{H})\leq \sum_{i=1}^{k}\gamma(\mathcal{H}_{i})=2\sum_{i=1}^{k}\nu(\mathcal{H}_{i})\le2\nu(\mathcal{H}).$$
To obtain the required equality, it suffices to show  that  $\gamma(\mathcal{H})\ge \sum_{i=1}^{k}\gamma(\mathcal{H}_{i})$ and
$\nu(\mathcal{H})\le \sum_{i=1}^{k}\nu(\mathcal{H}_{i})$.
Let $D$ be  a minimum dominating set of $\mathcal{H}$  such that $D$ contains vertices in $X$ as many as possible.
For partial hypergraph $\mathcal{H}_i$, by the construction of $\mathcal{\widehat{H}}_3$,
 there exists a degree-$1$ vertex $x_{i_si_t}^1$ in $H_i$ for any $x_{i_s}, x_{i_t}\in X_i$.  Because, by assumption, $D$
 contains vertices in $X$ as many as possible, this implies that $|D\cap X_i|\geq |X_i|-1$.
  On the other hand, by Proposition \ref{pro3.0}, we have $\gamma(\mathcal{H}_i)=|X_i|-1$.
So $\gamma(\mathcal{H})=|D|\ge\sum_{i=1}^k|D\cap X_i|\geq \sum_{i=1}^{k}\gamma(\mathcal{H}_{i})$.
Similarly, let $M$ be a maximum matching of $\mathcal{H}$ such that $M$ contains the edges in
$\cup_{i=1}^{k}\mathcal{H}_{i}$ as many as possible. We claim that $M\subseteq \cup_{i=1}^{k}\mathcal{H}_{i}$. Suppose not, there exists an edge  $e\in (\mathcal{Z}_1\cup \mathcal{Z}_2)\cap M.$ Without loss of generality, assume  $e=\{x_{m_1},x_{m_2},x_n\}\in \mathcal{Z}_1$, where $x_{m_1},x_{m_2},\in X_m, x_n\in X_n$.
Replacing edge $e$ by the edge $\{x_{m_1},x_{m_1m_2}^1,x_{m_2}\}$ containing $x_{m_1}$ and $x_{m_2}$ in $\mathcal{H}_m$, we obtain a maximum matching
$M'$ of $\mathcal{H}$ that contains more edges in $\cup_{i=1}^{k}\mathcal{H}_{i}$  than $M$, a contradiction. Thus $\sum_{i=1}^{k}\nu(\mathcal{H}_{i})\geq |M|=\nu(\mathcal{H})$.
This establishes the equality $\gamma(\mathcal{H})=2\nu(\mathcal{H})$. Therefore, $\mathcal{H}\in \mathfrak{H}_3$.

Conversely, suppose that $\mathcal{H}\in \mathfrak{H}_3$, i.e., $\gamma(\mathcal{H})=2\nu(\mathcal{H})$, we show that
 $\mathcal{H}\in \mathcal{G}_3$.
 By Lemma \ref{lem3.1},
there exists a spanning partial
hypergraph $\mathcal{H}^{*}$ of $\mathcal{H}$ such that every edge in $\mathcal{H}^{*}$ contains
exactly one degree-$1$ vertex and
  $\gamma(\mathcal{H}^{*}_{i})=2\nu(\mathcal{H}^{*}_{i})$
for each component $\mathcal{H}^{*}_{i}$ of $\mathcal{H}^{*}$.
Let $\mathcal{H}^*_{1}, \mathcal{H}^*_{2}, \ldots, \mathcal{H}^*_{k}$ be all the connected  components of $\mathcal{H}^*$ and let $\mathcal{H}_i=\mathcal{H}[V(\mathcal{H}^{*}_{i})]$.
By Lemma \ref{thm3.1}, $\mathcal{H}_i=\mathcal{F}$ or $\mathcal{H}_i\in \mathcal{\widehat{H}}_{3}$.
Let $V(\mathcal{H}_i)=X_i\cup X_i(A_i)$ for all $1\le i\le k$.
To complete the proof of necessity,
it suffices to show that $\mathcal{H}$ is obtained from $\cup_{i=1}^k\mathcal{H}_i$ by adding some edges of $\mathcal{Z}_1\cup \mathcal{Z}_2$.
In other word,
  we show that $e\in \mathcal{Z}_1\cup \mathcal{Z}_2$ for each $e\in E(\mathcal{H})\setminus \cup_{i=1}^{k} E(\mathcal{H}_{i})$.
For each $1\le i\le k$, let $M_i$ be a maximum matching of $\mathcal{H}_i$.

 We claim that
there exists an $\mathcal{H}_{s}$, $1\leq s\leq k$, such that $|e\cap V(\mathcal{H}_{s})|=2$. Suppose not,
then $|e\cap V(\mathcal{H}_{i})|\leq1$ for all $1\leq i\leq k$.
Without loss of generality, we may assume that $e\cap V(\mathcal{H}_{1})=\{u_1\}$, $e\cap V(\mathcal{H}_{2})=\{u_2\}$ and
 $e\cap V(\mathcal{H}_{3})=\{u_3\}$.
 By Proposition \ref{pro3.1}, there exists a maximum matching $M_{i}$ of $\mathcal{H}_{i}$ such that
 $V(M_{i})\subseteq V(\mathcal{H}_{i})\setminus \{u_i\}$ for each $1\le i\le 3$.
 By Lemma \ref{lem3.3}, $\nu(\mathcal{H})=\sum_{i=1}^{k}\nu(\mathcal{H}_i)$.  Then $\nu(\mathcal{H})=\sum_{i=1}^k|M_i|$. Clearly,
 $\cup_{i=1}^{k} M_{i}\cup \{e\}$
 is a matching of $\mathcal{H}$, this is a contradiction.
 Without loss of generality, we may assume that
 $e\cap V(\mathcal{H}_1)=\{u_1, u_2\}$ and $e\cap V(\mathcal{H}_2)=\{u_3\}$.

 We further claim that $e\cap X_1=\{u_1, u_2\}$. Suppose not, then $|e\cap X_1|\le 1$.
 If $e\cap X_{1}=\emptyset$,
then, by Proposition \ref{pro3.1}, there exists a maximum matching $M_{1}$ of $\mathcal{H}_{1}$ such that
 $u_1,v_1\not\in V(M_{1})$,
 and a maximum matching $M_{2}$ of $\mathcal{H}_{2}$
 such that $u_3\not\in V(M_{2})$.
 Then $\cup_{i=1}^{k}M_i\cup \{e\}$
 is matching of $\mathcal{H}$.
 This contradicts the fact that $\nu(\mathcal{H})=\sum_{i=1}^{k}\nu(\mathcal{H}_i)$ in Lemma \ref{lem3.3}.
 So we may assume $|e\cap X_1|=1$ and let $u_1\in X_1$. Since each $|X_i|$ is odd, $|X_1|=3$
or $|X_1|\ge 5$.
  If $|X_1|\geq 5$ or $|X_1|=3$ and $u_1$ and $u_2$ are adjacent, then by Proposition \ref{pro3.1}, there exists a maximum matching $M_{1}$ of $\mathcal{H}_{1}$ such that
 $u_1, u_2\not\in V(M_{1})$, and a maximum matching $M_{2}$ of $\mathcal{H}_{2}$
 such that $u_3\not\in V(M_{2})$.
 Then $\cup_{i=1}^{k} M_{i}\cup\{e\}$
 is a matching of $\mathcal{H}$, contradicting Lemma \ref{lem3.3} again. Finally, we may assume that
 $|X_1|=3$ and $u_1$ and $u_2$ are not adjacent. Let $X_1\backslash \{u_1\}=\{x_{s}, x_{t}\}$. Then $u_2\in X_{st}$.
 If $|X_{st}|\ge 2$, then $\mathcal{H}_1\neq \mathcal{F}$, so $\mathcal{H}_1\in \mathcal{\widehat{H}}_3$.
Let $x_{st}^p\in X_{st}\setminus \{u_2\}$ where $p\ge 1$.
 Then we can choose a maximum matching $M_{1}=\{x_{s},x_{st}^p, x_{t}\}$ of $\mathcal{H}_{1}$ containing no
  $u_2$. By Proposition \ref{pro3.1}, there exists a maximum matching $M_{2}$ of $\mathcal{H}_{2}$
 such that $u_3\not\in V(M_{2})$.
 Then $\cup_{i=1}^{k} M_{i}\cup\{e\}$
 is a matching of $\mathcal{H}$, a contradiction. Thus,
 we have $|X_{st}|=1$.
  Choose a vertex $v_i\in X_i$ for  $2\leq i\leq k$, let $D=\bigcup_{i=2}^k(X_i\setminus \{v_i\})\cup\{u_1\}$.
It is easy to see that $D$ is a dominating set of $\mathcal{H}$.
By Proposition \ref{pro3.0} and Lemma \ref{lem3.3}, $\gamma(\mathcal{H})\leq
\sum_{i=1}^k(|X_i|-1)-1=\sum_{i=1}^k\gamma(\mathcal{H}_i)-1\le 2\nu(\mathcal{H})-1$, contradicting the assumption
that $\mathcal{H}\in \mathfrak{H}_3$. Consequently, $|e\cap X_1|=2$.

If $u_3\in X_2$, then $e\in \mathcal{Z}_1$, as desired. Otherwise, $u_3\in X_2(A_2)$. We show that $u_3\in S_2$.
Suppose not, let $u_3=x_{s_1s_2}^1$. Then, for all $x_{s_1s_2}^p\in X_{s_1s_2}$, we have $d_{\mathcal{H}_2}(x_{s_1s_2}^p)\ge 2$,
since $\mathcal{H}_2\in \mathcal{\widehat{H}}_3\cup \mathcal{F}$.
Choose a vertex $v_i\in X_i, 1\leq i\leq k, i\neq 2$.
Then $D=\bigcup_{i=1, i\neq 2}^k(X_i\setminus \{v_i\})\cup(X_2\setminus \{x_{s_1}, x_{s_2}\})$ is a dominating set of $\mathcal{H}$. As noted above,
we would obtain that $\gamma(\mathcal{H})\leq 2\nu(\mathcal{H})-1$, a contradiction.
Thus $u_3\in S_2$, and hence $e\in \mathcal{Z}_2$. The necessity follows.
This completes the proof of Theorem \ref{thm3.7}.
\qed




\begin{thebibliography}{10}

\bibitem{Acharya1} B.D. Acharya, Domination in hypergraphs, AKCE J. Combin. 4 (2007) 117--126.


\bibitem{Acharya2} B.D. Acharya, Domination in hypergraphs: II ¨C new directions, in: Proc. Int. Conf. ¨C ICDM, 2008, pp. 1--16.

\bibitem{RON}  R. Aharoni, Ryser's conjecture for tri-partite $3$-graphs, Combinatorica 21 (2001) 1--4.

\bibitem{Alon} N. Alon,  Transversal numbers of uniform hypergraphs, Graphs Combin. 6 (1990) 1--4.

\bibitem{Aru} S. Arumugam, B. Jose, C. Bujt\'{a}s, Zs. Tuza, Equality of domination and transversal numbers
in hypergraphs, Discrete Appl. Math. 161 (2013) 1859--1867.

\bibitem{Bujt} Cs. Bujt\'{a}s, M.A. Henning, Zs. Tuza, Transversals and domination in uniform hypergraphs, European J. Combin. 33 (2012) 62--71.


\bibitem{Coc} E.J. Cockayne, S.T. Hedetniemi, P.J. Slater, Matchings and transversals in hypergraphs, domination and independence-in trees, J. Combin. Theory Ser. B 27 (1979) 78--80.


\bibitem{Chv} V. Chv\'{a}tal, C. McDiarmid, Small transversals in hypergraphs, Combinatorica 12 (1992) 19--26.

\bibitem{Dorf} M. Dorfling, M.A. Henning, Linear hypergraphs with large transversal number
and maximum degree two, European J. Combin. 36 (2014) 231--236.

\bibitem{fu}Z. F\"uredi, Matchings and covers in hypergraphs, Graphs  Combin. 4 (1988) 115--206.



\bibitem{has} P. Haxell, A. Scott, On Ryser's Conjecture, Electron. J. Combin. 19 (2012) \#P23.

\bibitem{hhs1} T.W. Haynes, S.T. Hedetniemi, P.J. Slater (Eds.), Fundamentals of Domination in Graphs, Marcel Dekker, Inc., New York, 1998.

\bibitem{hhs2}T.W. Haynes, S.T. Hedetniemi, P.J. Slater (Eds.), Domination in Graphs: Advanced Topics, Marcel Dekker, Inc., New York, 1998.

\bibitem{Henning1} M.A. Henning, C. L\"{o}wenstein, Hypergraphs with large domination number
and edge sizes at least 3, Discrete Appl. Math. 160 (2012) 1757--1765.

\bibitem{Henning2} M.A. Henning, C. L\"{o}wenstein, Hypergraphs with large transversal number
and with edge sizes at least four, Cent. Eur. J. Math.  10(3) (2012) 1133--1140.


\bibitem{Henning3} M.A. Henning, A. Yeo, Hypergraphs with large transversal number and with
edge sizes at least three, J. Graph Theory 59 (2008) 326--348.

\bibitem{Henning4} M.A. Henning, A. Yeo, Transversals and matchings in $3$-uniform hypergraphs,
European J. Combin. 34 (2013) 217--228.



\bibitem{hy}M.A. Henning, A. Yeo, Total Domination in Graphs, Springer, New York, 2013.

\bibitem{Jose} B.K. Jose, Zs. Tuza, Hypergraph domination and strong independence,  Appl. Anal. Discrete Math. 3 (2009) 237--358.

\bibitem{Kano} M. Kano, Y. Wu, Q. Yu, Star-uniform graphs, Graphs  Combin.  26 (2010) 383--394.

\bibitem{kld}L. Kang, S. Li, Y. Dong, E. Shan, Matching and domination numbers in $r$-uniform
hypergraphs, J. Comb. Optim. DOI: 10.1007/s10878-016-0098-5

\bibitem{Lai} F.C. Lai, G. J. Chang, An upper bound for the transversal numbers of
$4$-uniform hypergraphs, J. Combin. Theory Ser. B 50 (1990) 129--133.











\end{thebibliography}
\end{document}